\documentclass[a4paper, 12pt, pointlessnumbers]{article}

\usepackage{graphicx}
\usepackage{mathrsfs}

\usepackage[T1]{fontenc}							
\usepackage[latin1]{inputenc}	

\usepackage[all]{xy}
\usepackage{pb-diagram, pb-xy}

\usepackage{amssymb}
\usepackage{ntheorem} 

\theoremstyle{plain}
\theoremseparator{.}

\newtheorem{defn}{\hspace{\parindent}Definition}
\newtheorem{cor}[defn]{\hspace{\parindent}Corollary}

\newtheorem{lem}[defn]{\hspace{\parindent}Lemma}
\newtheorem{thm}[defn]{\hspace{\parindent}Theorem}

\newcommand{\medrowheight}{\rule[-0.4em]{0em}{1.4em}}

\newcommand{\bbN}{{\mathbb N}}

\newcommand{\frakg}{{\mathfrak g}}
\newcommand{\proof}{{\em Proof. }}

\begin{document} 

\title{\normalfont \large \bfseries
On the tensor square of irreducible representations of reductive
Lie superalgebras}
\author{\normalsize T. Kr\"amer and R. Weissauer}
\date{}
\maketitle

The algebraic super representations of a semisimple Lie superalgebra~$\frakg$ over an algebraically closed field $k$ of $\mathit{char}(k)=0$ form a semisimple tensor category if and only if the Lie superalgebra is a product of simple classical Lie algebras and orthosymplectic Lie superalgebras~$\mathfrak{osp}(1|2m)$, see \cite{DH}. For such Lie superalgebras~$\frakg$ we classify all irreducible nontrivial super representations $V$ whose 
symmetric square $S^2(V)$ or alternating square $\Lambda^2(V)$
is irreducible or decomposes into an irreducible and a trivial
representation, a situation that naturally arises in~\cite{KW}, \cite{We1}, \cite{We2}.  Obviously we may assume~$\frakg$ to be simple. For the classification  we make use of a list of small representations from~\cite{AEV}, extended  to the super case. 

\bigskip

For a uniform treatment including the supercase~$\frakg=\mathfrak{osp}(1|2m)$ with the root system~$BC_m$, the terms representation, dimension, symmetric square, trace etc.~will always be understood in the super sense. 
For 
$A_m$, $B_m$, $C_m$, $BC_m$ with $m\geq 1$ and for $D_m $ with $m\geq 3$, let $st$ denote the standard representation
of dimension $m+1$, $2m+1$, $2m$, $2m-1$ and $2m$ respectively. With respect to the dominant fundamental weights $\beta_1, \dots, \beta_m$ for $\frakg$ as introduced below the representation $st$ has highest weight $\beta_1$, respectively $2\beta_1$ in the case of $B_1$.

\bigskip

In the formulation of theorem~\ref{thm:main_thm} we use exceptional isomorphisms. 
Thus the representation of $A_1=C_1$ of highest weight $2\beta_1$ is the representation~$st$ of~$B_1$, and the one of $B_1$ of highest weight $\beta_1$ is the representation $st$ of $A_1=C_1$.
For $B_2=C_2$ the second fundamental representation of the one is the standard representation of the other, and ditto for $A_3=D_3$. 

\bigskip

We will assume~$dim(V)\geq 0$ since for the parity flip  $\Pi$ of the super structure   we have $dim(\Pi(V))=-dim(V)$ and 
$S^2(\Pi(V))\cong \Lambda^2(V)$.

\begin{thm} \label{thm:main_thm}
Suppose that $V$ is a nontrivial irreducible representation of~$\frakg$ with $dim(V)\geq 0$, and let $\lambda$ be its highest weight. 
\begin{enumerate}
 \item
If $S^2(V)$ is irreducible, then up to duality one has $V=st$, and $\frakg$ is of type $A_m$, $C_m$ or $BC_m$.
If $S^2(V)$ splits as a trivial plus a nontrivial irreducible representation, then $V=st$ and $\frakg$ is of type $B_m$ or $D_m$, or up to duality one of the following cases occurs 
\[
 \begin{array}{|c||c|c|c|c|} \hline
  \frakg & B_3 & D_4 & E_6 & G_2 \\ \hline
  \lambda & \beta_3 & \beta_3, \beta_4 & \beta_1 & \beta_1 \\ \hline
 \end{array}
\]
\item
If $\Lambda^2(V)$ is irreducible, then up to duality one has~$V=st$, and $\frakg$ is of type $A_m$, $B_m$ or $D_m$,
or up to duality one of the following cases occurs
\[
 \begin{array}{|c||c|c|c|c|c|} \hline
  \frakg & \multicolumn{2}{c|}{A_m} & D_4 & D_5 & E_6 \\ \hline
  \lambda & 2\beta_1 & \beta_2, \, m\geq 2 & \beta_3, \beta_4 & \beta_4, \beta_5 & \beta_1 \\ \hline
 \end{array}
\]
If $\Lambda^2(V)$ splits as a trivial plus a nontrivial irreducible representation, then $V=st$ and $\frakg$ is of type $C_m \, (m\geq 2)$ or $BC_m$, or one of the following cases occurs
\[
 \begin{array}{|c||c|c|c|c|c|} \hline
  \frakg & A_1 & A_5 & C_3 & D_6 & E_7 \\ \hline
  \lambda & 3\beta_1 & \beta_3 & \beta_3 & \beta_5, \beta_6 & \beta_7 \\ \hline
 \end{array}
\]
\end{enumerate}
\end{thm}

\bigskip

{\em Example 1.} 
If $\frakg$ is of type $A_1$, the irreducible representations of $\frakg = \mathfrak{sl}(2)$ are of the form $V=S^n(st)$ for a unique weight $n\in \bbN_0$. In this case the symmetric square
 $S^2(V)=\bigoplus_{i=0,1,\dots,\lfloor n/2\rfloor} S^{2n-4i}(st)$ is irreducible for $n=1$, irreducible up to a trivial representation for $n=2$, and neither of these otherwise. 
The alternating square $\Lambda^2(V)=\bigoplus_{i=0,1,\dots, \lfloor (n-1)/2 \rfloor} S^{2(n-1)-4i}(st)$ is irreducible for~$n\in \{1,2\}$, irreducible up to a trivial representation for~$n=3$, and neither of these otherwise.

\bigskip

{\em Example 2.} If $\frakg$ is of type $BC_1$, by~\cite{Dj2} for any irreducible representation~$V$ of $\frakg=\mathfrak{osp}(1|2)$ with $\dim(V)>0$ there exists a unique weight $n\in \bbN_0$ such that the even resp.~odd parts of $V=V_0 \oplus V_1$, as representations of~$\frakg_0 = \mathfrak{sl}(2)$, are $V_0 = S^n(st)$	and $V_1 = S^{n-1}(st)$. The even part of $S^2(V)$, as a representation of $\mathfrak{sl}(2)$, is given by
\[
	(S^2(V))_0 \;=\; S^2(S^n(st)) \oplus \Lambda^2(S^{n-1}(st))
	\;=\; S^{2n}(st) \; \oplus \; 2 \cdot \bigoplus_{i=1}^{\lfloor \frac{n}{2} \rfloor} S^{2n-4i}(st),
\]
%
so $S^2(V)$ is irreducible for $n=1$, but for $n>1$ it contains at least three irreducible summands.
Similarly
\[
	(\Lambda^2(V))_0 \;=\; \Lambda^2(S^n(st)) \oplus S^2(S^{n-1}(st))
	\;=\; 2 \cdot \bigoplus_{i=0}^{\lfloor \frac{n-1}{2} \rfloor} S^{2(n-1)-4i}(st),
\]
so $\Lambda^2(V)$ contains at least two nontrivial irreducible summands for $n>1$. For $n=1$ it is obvious that $\Lambda^2(V)$ splits into the trivial and a non-trivial irreducible representation. 

\bigskip

In view of the above two examples, in what follows we may and will  assume that $\frakg$ is neither of type $A_1$ nor of type $BC_1$.

\pagebreak

{\em Some notations.} There exists an invariant nondegenerate bilinear form~$(\cdot, \cdot)$ on $\frakg$. Any other invariant form on $\frakg$ is a multiple of this one (see~\cite{Sc} p.~94 for the case of~$\mathfrak{osp}(1|2m)$). 
Let~$\alpha_1, \dots, \alpha_m$ be a system of simple roots for~$\frakg$, and let $\beta_1, \dots, \beta_m$ denote the fundamental dominant weights given by $( \alpha_i^\vee, \beta_j) = \delta_{ij}$ for $\alpha_j^\vee = 2\alpha_j/( \alpha_j, \alpha_j)$. Then $\rho = \beta_1+\cdots +\beta_m$ is the half-sum, resp.~half-super-sum as in~\cite{CT} in the case of $\mathfrak{osp}(1|2m)$, of all positive roots. 

\begin{lem} \label{lem:ineq}
For simple $\frakg \neq A_1, BC_1$, and with the index set $I$ as given in table~\ref{tab:numerical}, one has
$$ 2\, \vert \rho\vert \, max_{i\in I}(\vert \alpha_i^\vee\vert) \;<\; dim(\frakg) - 1.$$
\end{lem}

\proof Notice that $|\rho|$ as well as the $|\alpha_i^\vee|$ depend on the chosen form $(\cdot, \cdot)$. However, the product $2\, |\rho| \, max_{i\in I} (|\alpha_i^\vee|)$ only depends on the root system. To compute it choose an isometric embedding of the dual of a Cartan algebra of $\frakg$ into Euclidean space 
(with standard basis $\epsilon_1, \epsilon_2, \dots$) as in~\cite{Bo}. For $\mathfrak{osp}(1|2m)$ we choose a similar embedding following~\cite{CT}, p.~353f. The lemma then is an immediate consequence of 
table~\ref{tab:numerical}.  
\\

\begin{table}[h] 
\footnotesize
\[
\begin{array}{|c||c|c|c|} \hline \medrowheight
	& |\rho|^2 & max_{i\in I} |\alpha_i^\vee| & \dim(\frakg) 
	\\ \hline \medrowheight
	A_m & \frac{1}{12} m(m+1)(m+2) & \sqrt{2} & m(m+2)
	\\ \hline \medrowheight
	B_m & \frac{1}{12} m(2m-1)(2m+1) & 2 & m(2m+1)
	\\ \hline \medrowheight
	C_m & \frac{1}{6} m(m+1)(2m+1) & \sqrt{2} & m(2m+1)
	\\ \hline \medrowheight
	D_m & \frac{1}{6} m(m-1)(2m-1)& \sqrt{2} & m(2m-1)
	\\ \hline \medrowheight
	E_6 & 78 & \sqrt{2} & 78\\ \hline \medrowheight
	E_7 & \frac{399}{2} & \sqrt{2} & 133\\ \hline \medrowheight
	E_8 & 620 & \sqrt{2} & 248\\ \hline \medrowheight
	F_4 & 39 & 2 & 52\\ \hline \medrowheight
	G_2 & 14 & \sqrt{2} & 14\\ \hline \medrowheight
	BC_m & \frac{1}{12} m(2m-1)(2m+1) & \sqrt{2} & m(2m-1) \\ \hline 
\end{array}
\]
	\caption{Some numerical values depending on the type of $\frakg$. Concerning the middle column, $I := \{1,2,\dots, m\}$ in the non-super case, but in the case of $\mathfrak{osp}(1|2m)$ we here only consider the~{\em even} coroots.}
\label{tab:numerical}
\end{table}

{\em A list of small representations.} We need to extend the classification of small representations in~\cite{AEV} to include the super case. Notice that in the case of~$\frakg = BC_1$ one has~$\dim(V)=\dim(\frakg)=1$ for
all irreducible representations~$V$ with $\dim(V)\geq 0$.

\begin{lem} \label{lem:small_reps} For simple ~$\frakg \neq BC_1$ 
one has $0\leq \dim(V_\lambda)\leq \dim(\frakg)$ iff $\lambda$ appears in tables~\ref{tab:replist1}, \ref{tab:replist2} or~\ref{tab:replist3} below.
\end{lem}

\proof See~\cite{AEV} for the case of ordinary simple Lie algebras. 
For the case~$\frakg = \mathfrak{osp}(1|2m)$ with $m\geq 2$ we use the Kac-Weyl superdimension formula as given in~\cite{CT}.
Following loc.~cit.~p.~356, the irreducible representations of~$\mathfrak{osp}(1|2m)$ of nonnegative dimension are parametrized by those dominant weights which, with respect to the chosen embedding in Euclidean space, 
take the form~$\lambda = (\lambda_1, \dots, \lambda_m)$ with~$\lambda_1, \dots, \lambda_{m}\in \bbN_0$ and $\lambda_1\geq \cdots \geq \lambda_m$~. Notice that these are not all the dominant weights but only those in which the last fundamental weight $\beta_m = (\frac{1}{2}, \dots, \frac{1}{2})$ enters with {\em even} multiplicity. For such~$\lambda$, by the Kac-Weyl formula, the superdimension of the corresponding irreducible representation $V_\lambda$ is
\[
 \dim(V_\lambda) \; \; \;=\;
 \prod_{1\leq i<j\leq m} \Bigl( \frac{\lambda_i - \lambda_j}{j-i}+1\Bigr)
 \; \; \cdot
 \prod_{1\leq i<j\leq m}
 \Bigl( \frac{\lambda_i + \lambda_j}{2m+1-i-j} + 1 \Bigr).
\]

If $\lambda_1\geq 2$, the second product is $\geq 2$. So the classical Weyl formula applied to the first product shows that $\dim(V_\lambda)$ is at least twice the dimension of the simple $\mathfrak{sl}(m)$-module of highest weight $\lambda' = (\lambda_1-\lambda_m, \dots, \lambda_{m-1}-\lambda_m)$. Since $\dim(\mathfrak{sl}(m)) \geq 2\dim(\mathfrak{osp}(1|2m))$, we are then reduced to $\lambda'$ being in the list for $A_{m-1}$ in table~\ref{tab:replist1}.
But $\lambda$ is obtained from $\lambda'$ by adding a multiple of~$2\beta_m = \epsilon_1+\cdots + \epsilon_m$. Replacing $\lambda$ by $\lambda + 2\beta_m$  in the Kac-Weyl formula leaves the first product unchanged but strictly increases the second one. So there only remain finitely many cases which are in table~\ref{tab:replist3}.

If $\lambda_1 = 1$, then $\lambda = \epsilon_1 + \cdots + \epsilon_r$ for some $r\leq m$. Then by~\cite{Dj2} the even resp.~odd parts of $V_\lambda$ are $\Lambda^r(k^{2m})$ resp.~$\Lambda^{r-1}(k^{2m})$, so $\dim(V_\lambda) = {2m\choose r} - {2m\choose r-1}$, and from this one can deduce that all relevant cases are in table~\ref{tab:replist3}.
\\

\begin{cor} \label{cor:dimension_2} 
Let $V$ be a nontrivial irreducible representation of $\frakg$ with
$dim(V)\geq 0$, and assume $\frakg \neq A_1, BC_1$.  Then 
\[ \dim(V)\geq 2, \]
with equality holding only in case  $\frakg= \mathfrak{osp}(1|4)$  and $V=V_{\epsilon_1+\epsilon_2}$.
\end{cor}

\bigskip

{\em The index.} Assume $\frakg$ to be simple.
For any representation $\varphi: \frakg\rightarrow \mathfrak{gl}(V)$ 
the form 
	$\mathrm{tr}(\varphi(X)\circ \varphi(Y))$ is a nondegenerate invariant form on $\frakg$,
	hence equal to $ l(\varphi)\cdot ( X, Y)$ for a unique constant $l(\varphi)$,
called the {\em index} of $\varphi$. Obviously the index is additive under direct sums and functorial in the sense that for 
any representation $\psi: \mathfrak{gl}(V) \to \mathfrak{gl}(V')$ there exists
a constant $l(\psi,V)$ such that $l(\psi\circ\varphi) = l(\psi,V)\cdot l(\varphi)$. If $\psi$ is the
second symmetric or alternating power, this becomes
\[ l(S^2, V) = \dim(V) + 2
 \qquad
 \textnormal{and}
 \qquad
   l(\Lambda^2, V) = \dim(V) - 2 \, . \]

\pagebreak
To see this for $\psi = S^2$, write $V=k^{r|s}$ and notice $\varphi(\frakg)\subset \mathfrak{sl}(r|s)$ since~$\frakg$ is simple. Now~$r\neq s$, so $\mathfrak{sl}(r|s)$ is also simple by~\cite{Sc} p.~128. Thus it \mbox{suffices} to check for a suitable elementary matrix $X = e_{ii}-e_{jj} \in \mathfrak{sl}(r|s)$ that $\mathrm{tr}((S^2(X))^2)/\mathrm{tr}(X^2) = r-s+2$. The case where~$\psi = \Lambda^2$ may be treated similarly or reduced to the case $\psi = S^2$ by a parity flip.

\bigskip

 The Casimir operator acts by a scalar on any irreducible representation~$V_\mu$
 with highest weight $\mu$. This scalar is proportional to
 $$	 
	c(\mu) \; =\; ( \mu, \mu )
	+ 2( \mu, \rho ) \, ,
\ $$
as may be seen by  looking at the action on a highest weight vector (this works for $\mathfrak{osp}(1|2m)$ just as in the classical case; see the top of p.~28 in~\cite{Dj2} and p.~223 in~\cite{Dj1} for the relevant setting). Notice $c(\mu)>0$ for dominant weights $\mu\neq 0$. Hence from the definition of the index we obtain for some universal constant $\kappa := \dim(Ad)\cdot c(Ad) \neq 0$, depending only on  $\frakg$ and on $(\cdot, \cdot )$, that
$$
\kappa \cdot l(V_\mu) \;=\;  \dim(V_\mu) \cdot c(\mu) \, .$$

\bigskip

{\em Proof of theorem~\ref{thm:main_thm} (the symmetric case)}.
By assumption $\dim(V)\geq 0$, so $V=V_\lambda$ is a highest weight representation for some dominant weight $\lambda$.
By corollary \ref{cor:dimension_2}
we may also assume $n := \dim(V_\lambda) > 2$.
Suppose that
$S^2(V_\lambda)$ is irreducible or splits into a (then one-dimensional) trivial plus one further irreducible summand; put~$\delta = 0$ resp.~$\delta = 1$ in these two cases. Then we claim
$$
	(\dim(V_\lambda) -2\delta)\cdot \vert \lambda \vert^2 \;=\; 2\,( \lambda, \rho) \, .
$$
Since $S^2(V_\lambda)=V_{2\lambda}\oplus k^{\delta}$ and hence $\dim(V_{2\lambda})=
n(n+1)/2 - \delta$, this follows 
from the equality $(n+2)l(V_\lambda) =l(V_{2\lambda})$ for $l(S^2(V_\lambda))=l(V_{2\lambda})$.
Indeed, via the above Casimir formula for the index then $(n+2)nc(\lambda)=(n(n+1)/2-\delta)c(2\lambda)$, hence $(n^2-4\delta)( \lambda, \lambda ) = 2(n+2\delta)( \lambda, \rho)$. So our claim follows from~$n+2\delta > 0$.

By the Cauchy-Schwartz inequality the claim above shows that
$$  |\lambda| \; \leq \; \frac{2\, |\rho|}{\dim(V_\lambda)-2\delta} \, ,$$
which by lemma \ref{lem:ineq}
implies 
$$ \vert(\lambda, \alpha_i^\vee ) \vert \ \leq\  |\lambda|  \cdot |\alpha_i^\vee|  
\ \leq \ \frac{ 2 \, \vert\rho\vert \, max_{i\in I}(\vert \alpha^\vee_i\vert)}{\dim(V_\lambda)-2\delta}  \ < \  \frac{ \dim( \frakg ) -1 }{\dim(V_\lambda)-2\delta} \,  .$$
Since for the dominant integral weight $\lambda\neq 0$ the coefficients $(\lambda, \alpha_i^\vee )$ are \mbox{integers} and not all of them are zero, it follows that $dim(V_\lambda)\leq\dim(\frakg )$. 
%
Hence the proof of our theorem for $S^2(V)$ is finished 
by tables~\ref{tab:replist1}, \ref{tab:replist2} and~\ref{tab:replist3}.

\bigskip

{\em Proof of theorem~\ref{thm:main_thm} (the alternating case)}. 
Now assume that
$\Lambda^2(V_\lambda)$ is \mbox{irreducible} or splits into a (then one-dimensional) trivial plus one further \mbox{irreducible} summand. Again put $\delta = 0$ resp.~$\delta = 1$ in these two cases. Write $\Lambda^2(V)=W\oplus k^\delta$.
Then $W$ is irreducible of dimension $n(n-1)/2 - \delta > 0$, hence we may write $W=V_\mu$ for some dominant weight $\mu$. 

For a highest weight vector $v$ of $V_\lambda$ there exists 
at least one simple positive root $\alpha=\alpha_i$ and $X_{-\alpha}\in \frakg_{-\alpha}$ such that
$w= X_{-\alpha}v \neq 0$, since otherwise $V_\lambda$ would be trivial. Now $v\wedge w\neq 0$ has the weight $2\lambda - \alpha$ which is non-zero for~$\frakg\neq A_1$, so $v\wedge w$ has to come from 
the non-trivial summand $W=V_\mu$. 

By the Poincare-Birkhoff-Witt theorem respectively its super-analogon, any weight of $V_\mu$ is in $\mu - {\cal P}$, where ${\cal P}$ is the monoid of all nonnegative integral linear combinations
of the positive simple roots $\alpha_1,..,\alpha_m$. Thus $2\lambda - \alpha =  \mu -  \alpha' $ for some $\alpha'\in {\cal P}$. 
As a weight 
of $V_\lambda\otimes V_\lambda$ also $\mu \in 2\lambda - {\cal P}$, 
because the weights of $V_\lambda$ are in $\lambda - {\cal P}$.
But $\mu\neq 2\lambda$, since $2\lambda$ has multiplicity one in $V_\lambda\otimes V_\lambda$ and is the highest weight of the symmetric square $S^2(V_\lambda)$.  
Hence $\mu\ =\ 2\lambda - \alpha'' $ for some $0\neq \alpha'' \in  {\cal P}$,
and by a comparison $\alpha = \alpha' + \alpha''$. Since $\alpha$ is a simple positive root this implies $\alpha'=0$, $\alpha=\alpha''$ and
$$ \mu \ =\ 2\lambda - \alpha \, .$$
Thus  $\alpha = \alpha_i$ is the unique (!)~simple root with $X_{-\alpha}v\neq 0$ for $X_{-\alpha}\in\frakg_{-\alpha}$.
For~$j\neq i$ then $h_{\alpha_j}v =[X_{\alpha_j},X_{-\alpha_j}]v=0$, so $(\lambda,\alpha_j^\vee)=0$. Hence
\[
\lambda = r\cdot \beta_i \, .
\]

On the other hand using $n>2$, a calculation similar to the one in the case of the symmetric square based on the relation $l(V_\mu)=l(\Lambda^2(V_\lambda)) = (n-2)\cdot l(V_\lambda)$ and the Casimir formula,
gives 
\[
 \Bigl(n-1-\frac{2\delta}{n}\Bigr)
 \Bigl( ( \alpha, \rho) + \frac{1}{2} \,|\alpha|^2 - |\lambda - \alpha|^2 \Bigr)
\;=\;
 \Bigl( 1 - \frac{2\delta}{n} \Bigr) \cdot c(\lambda) \, ,
\]
so that $2\delta/n < 1$ and $c(\lambda)>0$ imply $( \alpha, \rho) + \frac{1}{2}|\alpha|^2 - |\lambda - \alpha|^2 > 0$.  For simple positive roots $( \alpha, \rho) = \frac{1}{2}|\alpha|^2 $ and hence
$2(\alpha,\lambda) - \vert \lambda\vert^2\ > \ 0$.
For $\lambda = r\cdot \beta_i$ this implies
 $ 2(\alpha,\lambda) - \vert \lambda\vert^2 = r(\vert\alpha_i\vert^2  - r\vert\beta_i\vert^2 ) >0$, so the integer $r$ satisfies 
 $$ r \; < \; r_i \; := \; \frac{\vert \alpha_i\vert^2}{\vert \beta_i\vert^2} \, .$$
This only leaves the cases $1< r_i$ in table~\ref{tab:ri}. One then checks that  these are already listed in tables~\ref{tab:replist1} to~\ref{tab:replist3} below (for the case of $A_m$, one may here by duality of course assume $i<(m+1)/2$).
\\

\begin{table}[h] 
\small
\[
\begin{array}{|r||r|l|} \hline \medrowheight
	& & r_i
	\\ \hline  \medrowheight
	A_m & 1\leq i\leq m & \frac{2(m+1)}{i(m+1-i)} 
	\\ \hline \medrowheight
	B_m, C_m, BC_m & 1\leq i<m & \frac{2}{i} \\ \medrowheight
	&  i=m & \frac{4}{m}
	\\ \hline \medrowheight
	D_m & 1\leq i<m-1 & \frac{2}{i} \\ \medrowheight
	&  m-1\leq i\leq  m & \frac{8}{m}
	\\ \hline \medrowheight
	E_6 & i = 1,\dots, 6 & \frac{3}{2}, 1, \frac{3}{5}, \frac{1}{3}, \frac{3}{5}, \frac{3}{2} \\ \hline \medrowheight
	E_7 & i = 1, \dots, 7 & 1, \frac{4}{7}, \frac{1}{3},  \frac{1}{6},  \frac{4}{15}, \frac{1}{2}, \frac{4}{3} \\ \hline \medrowheight
	E_8 & i = 1, \dots, 8 & \frac{1}{2}, \frac{1}{4}, \frac{1}{7}, \frac{1}{15}, \frac{1}{10}, \frac{1}{6}, \frac{1}{3}, 1 \\ \hline \medrowheight
	F_4 & i = 1, \dots, 4 & 1, \frac{1}{3}, \frac{1}{3}, 1 \\ \hline \medrowheight
	G_2 & i = 1, 2 & 1 \\ \hline 
\end{array}
\]
\caption{The values of $r_i$ for the various types of $\frakg$.}
\label{tab:ri}
\end{table}

\vspace{6em}

\begin{table}[h]
\[
\begin{array}{|c||lr|r|r|} \hline  \medrowheight
  & & \lambda & S^2(V_\lambda) & \Lambda^2(V_\lambda) 
  \\ \hline 
  A_m & m = 1: & 2\beta_1 & V_{4\beta_1} \oplus k & V_{2\beta_1} \\
	& m\geq 2: & \beta_1 + \beta_m 
	& V_{2\beta_1 + 2\beta_m} \oplus V_{\beta_2 + \beta_{m-1}} 
	& V_{\beta_2 + 2\beta_m} \oplus V_{2\beta_1 + \beta_{m-1}}  \\
	&&& \oplus V_{\beta_1+\beta_m} \oplus k & \oplus V_{\beta_1 + \beta_m}
	\\ \hline 
B_m & m = 1: & 2\beta_1 & V_{4\beta_1} \oplus k & V_{2\beta_1} 
	\\ 
	& m = 2: & 2\beta_2 
	&  V_{\beta_1} \oplus V_{2\beta_1} \oplus V_{4\beta_2} \oplus k 
	&  V_{\beta_1+2\beta_2} \oplus V_{2\beta_2} \\
	& m =3: & \beta_2 
	& V_{2\beta_1} \oplus V_{2\beta_2} \oplus V_{2\beta_3} \oplus k & V_{\beta_1 + 2\beta_3} \oplus V_{\beta_2} \\
	& m \geq 4: & \beta_2
	& V_{2\beta_1} \oplus V_{2\beta_2} \oplus V_{\beta_4} \oplus k & V_{\beta_1 + \beta_3} \oplus V_{\beta_2}
	\\ \hline 
	C_m & m = 1: & 2\beta_1 & V_{4\beta_1} \oplus k & V_{2\beta_1} \\ 
	& m\geq 2: & 2\beta_1
	& V_{4\beta_1} \oplus V_{2\beta_2} \oplus V_{\beta_2} \oplus k & V_{2\beta_1} \oplus V_{2\beta_1 + \beta_2}
	\\ \hline 
	D_m & m = 3: & \beta_2 + \beta_3 
	& V_{2\beta_2 + 2\beta_3} \oplus V_{2\beta_1} 
	& V_{\beta_1 + 2\beta_3} \oplus V_{\beta_1+2\beta_2}  \\
	&&& \oplus V_{\beta_2+\beta_3} \oplus k & \oplus V_{\beta_2+\beta_3}\\
	& m = 4: & \beta_2
	& V_{2\beta_1} \oplus V_{2\beta_2} \oplus V_{2\beta_3} 
	& V_{\beta_2} \oplus V_{\beta_1+\beta_3+\beta_4} \\
	&&& \oplus V_{2\beta_4} \oplus k & \\
	& m \geq 5: & \beta_2
	& V_{2\beta_1} \oplus V_{2\beta_2} \oplus V_{\beta_4} \oplus k & V_{\beta_2} \oplus V_{\beta_1 + \beta_3}
	\\ \hline 
	E_6 && \beta_2 & V_{2\beta_2} \oplus V_{\beta_1+\beta_6} \oplus k & V_{\beta_2} \oplus V_{\beta_4}
	\\ \hline 
	E_7 && \beta_1 & V_{2\beta_1} \oplus V_{\beta_6} \oplus k & V_{\beta_1} \oplus V_{\beta3}
	\\ \hline 
	E_8 && \beta_8 & V_{\beta_1} \oplus V_{2\beta_8} \oplus k & V_{\beta_7} \oplus V_{\beta_8} 
	\\ \hline 
	F_4 && \beta_1 & V_{2\beta_1} \oplus V_{2\beta_4} \oplus k & V_{\beta_1} \oplus V_{\beta_2} 
	\\ \hline 
	G_2 && \beta_2 & V_{2\beta_1} \oplus V_{2\beta_2} \oplus k & V_{3\beta_1} \oplus V_{\beta_2} \\ \hline
\end{array}
\]
\caption{$S^2(V_\lambda)$ and $\Lambda^2(V_\lambda)$ in all classical cases where $\dim(V_\lambda) = \dim(\frakg)$; here the $V_\lambda$ are precisely the adjoint representations by~\cite{AEV}.}
\label{tab:replist1}
\end{table}

\begin{table}[h]
\[
\begin{array}{|c||lr|r|r|} \hline  \medrowheight
  & & \lambda & S^2(V_\lambda) & \Lambda^2(V_\lambda) 
  \\ \hline 
  A_m & m\geq 1: & \beta_1 & V_{2\beta_1} & V_{\beta_2}
	\\ 
	& & \beta_m & V_{2\beta_m} & V_{\beta_{m-1}}
	\\ 
	& m\geq 2:	& 2\beta_1 &  V_{4\beta_1}\oplus V_{2\beta_2} & V_{2\beta_1+\beta_2}
	\\ 
	& & 2\beta_m & V_{4\beta_m}\oplus V_{2\beta_{m-1}} & V_{2\beta_m + \beta_{m-1}}
	\\ 
	& m = 3: & \beta_2 & V_{2\beta_2}\oplus k & V_{\beta_1+\beta_3}
	\\ 
	& m \geq 4: & \beta_2 & V_{2\beta_2}\oplus V_{\beta_4} & V_{\beta_1 + \beta_3}
	\\ 
	&& \beta_{m-1} & V_{2\beta_{m-1}}\oplus V_{\beta_{m-3}} & V_{\beta_m+\beta_{m-2}}
	\\ 
	& m = 5: & \beta_3 & V_{2\beta_3}\oplus V_{\beta_1+\beta_5} & V_{\beta_2+\beta_4} \oplus k
	\\ 
	& m = 6,7: & \beta_3 & V_{2\beta_3}\oplus V_{\beta_1+\beta_5} & V_{\beta_2+\beta_4} \oplus V_{\beta_6}
	\\ 
	&&  \beta_{m-2} & V_{2\beta_{m-2}}\oplus V_{\beta_m+\beta_{m-4}} & V_{\beta_{m-1}+\beta_{m-3}} \oplus V_{\beta_{m-5}}
	\\ \hline 
B_m & m = 1: & \beta_1 & V_{2\beta_1} & k
	\\ 
	& m\geq 2: & \beta_1 & V_{2\beta_1}\oplus k & V_{\beta_2}
	\\ 
	& m = 2: & \beta_2 & V_{2\beta_2} & V_{\beta_1} \oplus k
	\\ 
	& m = 3: & \beta_3 & V_{2\beta_3}\oplus k & V_{\beta_1}\oplus V_{\beta_2}
	\\ 
	& m = 4: & \beta_4 & V_{2\beta_4}\oplus V_{\beta_1}\oplus k & V_{\beta_2}\oplus V_{\beta_3}
	\\ 
	& m=5: & \beta_m & V_{2\beta_m}\oplus V_{\beta_{m-3}} 
	\oplus V_{\beta_{m-4}} & V_{\beta_3} \oplus V_{\beta_4} \oplus k
	\\ 
	& m=6: & \beta_m & V_{2\beta_m}\oplus V_{\beta_{m-3}} 
	\oplus V_{\beta_{m-4}} & V_{\beta_1} \oplus V_{\beta_4} \oplus V_{\beta_5} \oplus k
	\\ \hline 
C_m & m = 1: & \beta_1 & V_{2\beta_1} & k \\
	& m\geq 2: & \beta_1 & V_{2\beta_1} & V_{\beta_2} \oplus k
	\\ 
	& m = 2: & \beta_2 & V_{2\beta_2} \oplus k & V_{2\beta_1}
	\\ 
	& m = 3: & \beta_2 & V_{\beta_2}\oplus V_{2\beta_2}\oplus k & V_{2\beta_1} \oplus V_{\beta_1+\beta_3}
	\\ 
	& m = 3: & \beta_3 & V_{2\beta_3} \oplus V_{2\beta_1} & V_{2\beta_2} \oplus k
	\\ 
	& m \geq 4: & \beta_2 & V_{\beta_4}\oplus V_{2\beta_2}\oplus V_{\beta_2} \oplus k & V_{\beta_1+\beta_3} \oplus V_{2\beta_1}
	\\ \hline 
	D_m & m\geq 3: & \beta_1 & V_{2\beta_1} \oplus k & V_{\beta_2}
	\\ 
	& m = 4: & \beta_3 & V_{2\beta_3} \oplus k & V_{\beta_2}
	\\ 
	&& \beta_4 & V_{2\beta_4} \oplus k & V_{\beta_2}
	\\ 
	& m = 5: & \beta_4 & V_{2\beta_{4}}\oplus V_{\beta_{1}} & V_{\beta_3}
	\\ 
	&& \beta_5 & V_{2\beta_{5}}\oplus V_{\beta_{1}} & V_{\beta_3}
	\\ 
	& m = 6: & \beta_5 & V_{2\beta_{5}}\oplus V_{\beta_{2}} & V_{\beta_4} \oplus k
	\\ 
	&& \beta_6 & V_{2\beta_{6}}\oplus V_{\beta_{2}} & V_{\beta_4} \oplus k
	\\ 
	& m = 7: & \beta_6 & V_{2\beta_{6}}\oplus V_{\beta_{4}} & V_{\beta_1} \oplus V_{\beta_5} \oplus k
	\\ 
	&& \beta_7 & V_{2\beta_{7}}\oplus V_{\beta_{4}} & V_{\beta_1} \oplus V_{\beta_5} \oplus k
	\\ \hline 
	E_6 && \beta_1 & V_{2\beta_1}\oplus k & V_{\beta_3}
	\\ 
	&& \beta_6 & V_{2\beta_6}\oplus k & V_{\beta_5}
	\\ \hline 
	E_7 && \beta_7 & V_{2\beta_7}\oplus V_{\beta_1} & V_{\beta_6} \oplus k
	\\ \hline 
	F_4 && \beta_4 & V_{2\beta_4}\oplus V_{\beta_4} \oplus k & V_{\beta_3} \oplus V_{\beta_1}
	\\ \hline 
	G_2 && \beta_1 & V_{2\beta_1}\oplus k & V_{\beta_1} \oplus V_{\beta_2} \\ \hline
\end{array}
\]
\caption{$S^2(V_\lambda)$ and $\Lambda^2(V_\lambda)$ in all classical cases where $\dim(V_\lambda) < \dim(\frakg)$.}
\label{tab:replist2}
\end{table}

\begin{table}[b]
\[
\begin{array}{|lr|r|r|} \hline \medrowheight
  \lambda & & S^2(V_\lambda) & \Lambda^2(V_\lambda) 
  \\ \hline 
	\mu_1 & m\geq 1: &  V_{2\mu_1} & V_{\mu_2} \oplus k
	\\ \hline 
	2\mu_1 & m =1: & V_{4\mu_1} \oplus  k \oplus \Pi V_{\mu_1}  & V_{2\mu_1} \oplus \Pi V_{3\mu_1}
	\\ \cline{3-4} 
	& m \geq 2: & V_{4\mu_1} \oplus V_{\mu_2} \oplus V_{2\mu_2} \oplus k  & V_{2\mu_1+\mu_2} \oplus V_{2\mu_1}
	\\ \hline 
	\mu_1 + \mu_2 & m = 2: & 2V_{2\mu_1 + \mu_2} \oplus V_{2\mu_1+2\mu_2}  
	& V_{4\mu_1} \oplus V_{3\mu_2} \oplus 2V_{\mu_2}
	\\ 
	&& \oplus 2V_{2\mu_1} \oplus \Pi V_{\mu_1 + 2\mu_2}  & \oplus V_{2\mu_1 + \mu_2} \oplus 2V_{2\mu_2} \oplus k\\
	&& \oplus \Pi V_{3\mu_1} \oplus 2\Pi V_{\mu_1 + \mu_2} & \oplus \Pi V_{3\mu_1 + \mu_2} \oplus \Pi V_{\mu_1+2\mu_2} \\
	&&& \oplus \Pi V_{\mu_1 + \mu_2} \oplus \Pi V_{\mu_1}\\ \hline
	\mu_2 & m = 2: & V_{\mu_2} \oplus V_{2\mu_2} \oplus k \oplus \Pi V_{\mu_1} & V_{2\mu_1} \oplus \Pi V_{\mu_1 + \mu_2}
	\\ \cline{3-4} 
	& m = 3: & V_{\mu_2} \oplus V_{2\mu_2} \oplus k \oplus \Pi V_{\mu_3} & V_{2\mu_1} \oplus V_{\mu_1 + \mu_3}
	\\ \cline{3-4} 
	& m \geq 4: & V_{\mu_2} \oplus V_{2\mu_2} \oplus V_{\mu_4} \oplus k & V_{2\mu_1} \oplus V_{\mu_1 + \mu_3}
	\\ \hline 
	\mu_3 & m = 3: &  V_{2\mu_1} \oplus V_{2\mu_3} &  V_{\mu_2} \oplus V_{2\mu_2} \oplus k\\
	&& \oplus V_{\mu_1 + \mu_3} \oplus \Pi V_{\mu_1 + \mu_2} & \oplus \Pi V_{\mu_1} \oplus \Pi V_{\mu_3} \oplus \Pi V_{\mu_2 + \mu_3}
	\\ \cline{3-4} 
	& m = 4: & V_{2\mu_1} \oplus V_{2\mu_3} \oplus V_{\mu_1+\mu_3}  
	& V_{\mu_2} \oplus V_{\mu_4} \oplus V_{\mu_2+\mu_4}
	\\ 
	&& \oplus \Pi V_{\mu_1+\mu_4} & \oplus V_{2\mu_2} \oplus k \oplus \Pi V_{\mu_3} \\ \hline
	\mu_4 & m = 4: & V_{2\mu_1} \oplus V_{\mu_4} \oplus V_{\mu_2 + \mu_4} 
	& V_{\mu_1+\mu_3} \oplus V_{2\mu_1} \oplus V_{2\mu_3}
	\\ 
	&& \oplus V_{\mu_2} \oplus V_{2\mu_4} \oplus k & \oplus \Pi V_{\mu_1 + \mu_4} \oplus \Pi V_{\mu_3 + \mu_4} 
	\\
	&& \oplus \Pi V_{\mu_3} \oplus \Pi V_{\mu_2+\mu_3} \oplus \Pi V_{\mu_1} & \oplus \Pi V_{\mu_1+\mu_2} \\ \hline 
	\mu_5 & m = 5: & V_{2\mu_1} \oplus V_{\mu_3+\mu_5} \oplus V_{2\mu_3} 
	& V_{\mu_2} \oplus V_{2\mu_4} \oplus V_{\mu_2+\mu_4} 
	\\ 
	&& \oplus V_{\mu_1+\mu_3} \oplus V_{\mu_1 + \mu_5} \oplus 2V_{\mu_5} 
	& \oplus V_{\mu_4} \oplus V_{2\mu_2} \oplus k  \\ 
	&& \oplus \Pi V_{\mu_1 + \mu_4} \oplus \Pi V_{\mu_1 + \mu_2}  
	& \oplus \Pi V_{\mu_2 + \mu_3} \oplus \Pi V_{\mu_3} \oplus \Pi V_{\mu_5} \\
	&& \oplus \Pi V_{\mu_3 + \mu_4}
	& \oplus \Pi V_{\mu_1} \oplus \Pi V_{\mu_4+\mu_5} \oplus \Pi V_{\mu_2 + \mu_5} \\ \hline	
\end{array}
\]
\caption{$S^2(V_\lambda)$ and $\Lambda^2(V_\lambda)$ for all representations $V_\lambda$ of $\frakg = \mathfrak{osp}(1|2m)$ such that $0\leq \dim(V_\lambda) \leq \dim(\frakg)$. To avoid extra case distinctions we put \mbox{$\mu_i = \epsilon_1 + \cdots + \epsilon_i$}, i.e.~$\mu_i = \beta_i$ for $i<m$ but $\mu_m = 2\beta_m$.
One has $\dim(V_\lambda) = \dim(\frakg)$ iff either $\lambda = 2\mu_1$ or $m=4$, $\lambda = \mu_3$. }
\label{tab:replist3}
\end{table}

\end{document}